\documentclass[12pt,reqno]{amsart}
\usepackage{amsfonts}

\newtheorem{theorem}{Theorem}
\theoremstyle{plain}

\newtheorem{case}{Case}

\newtheorem{lemma}{Lemma}
\newtheorem{notation}{Notation}

\newtheorem{remark}{Remark}

\numberwithin{equation}{section}
\input{tcilatex}
\begin{document}
\title[MHD instability]{Variational approach to nonlinear gravity-driven
instabilities in a MHD setting}
\author{Hyung Ju Hwang}
\address{Department of Mathematics, Duke University, Durham NC 27708, USA}
\email{hjhwang@math.duke.edu}
\urladdr{}
\date{}
\subjclass{}
\keywords{}
\dedicatory{}
\thanks{}

\begin{abstract}
We establish a variational framework for nonlinear instabilities in a
setting of the ideal magnetohydrodynamic (MHD) equations. We apply a
variational method to various kind of smooth steady states which are shown
to be nonlinearly unstable for both incompressible and compressible ideal
MHD equations. Destabilizing effect of compressibility is justified as well
as stabilizing effect of magnetic field lines arising in MHD dynamics, which
distinguishes from the Rayleigh-Taylor instability in the absence of
magnetic field lines.
\end{abstract}

\maketitle

\section{Introduction}

Magnetohydrodynamic equations (MHD) serve as an important model for fluid
and gas dynamics and hydromagnetic instability is a fundamental phenomenon
in nature, for instance, oceans, atmosphere, and plasma. MHD instability
possesses extensive applications in both laboratory plasmas and astrophysics
such as in nuclear fusions, compression of thin foils for X-ray production,
and stellar dynamics. Nevertheless, there have been not many analytical
results to date due to its structural complexity such as the presence of
shock waves. Furthermore, nothing has been known about nonlinear
instabilities for MHD equations despite its importance and variety of
instabilities.

The main purpose of this article is to present a variational framework in
the passage from linear to nonlinear instability in a setting of the ideal
MHD system and derive nonlinear instability around different steady states
for both incompressible and compressible ideal MHD equations.

We consider the equations of ideal magnetohydrodynamics(MHD) for inviscid
flows:

\begin{align}
\rho _{t}+\nabla \cdot \left( \rho V\right) & =0,  \notag \\
\rho \left( V_{t}+V\cdot \nabla V\right) & =\left( \nabla \times B\right)
\times B-\nabla P+\rho \vec{g},  \label{origin} \\
B_{t}& =\nabla \times \left( V\times B\right) ,  \notag \\
\nabla \cdot B& =0.  \notag
\end{align}%
Here $\rho $ is plasma density, $V\,$\ velocity, $B$ magnetic field, $P\,\ $%
plasma pressure, and $\vec{g}$ is the gravitational field. The $x$ axis is
taken along the gravitational field, which is assumed to be uniform:%
\begin{equation*}
\vec{g}=\left( g,0,0\right) .
\end{equation*}%
The condition for steady state $\left( \rho _{0},B_{0},P_{0}\right) $ with $%
V_{0}\equiv 0$ is then%
\begin{equation}
\frac{d}{dx}\left( P_{0}+\frac{1}{2}B_{0}^{2}\right) =\rho _{0}g.
\label{steady-state}
\end{equation}%
We assume for density profile that%
\begin{equation}
\min_{x}\rho _{0}\geq c>0.  \label{min}
\end{equation}%
Our domain is $t\geq 0$ and 
\begin{equation*}
D=\left\{ 0\leq x\leq 2\pi ,0\leq y\leq 2\pi \right\} .
\end{equation*}%
We assume periodic conditions at the boundary for $V$ and $B$. Both
incompressible and compressible fluids are considered: In the incompressible
case, 
\begin{equation*}
\nabla \cdot V=0.
\end{equation*}%
In the compressible case, $\gamma $ is an adiabatic index which relates the
pressure $P$ to the density $\rho $ by%
\begin{eqnarray}
P\left( \rho \right) &=&C\rho ^{\gamma },\   \label{p-rh} \\
\frac{\nabla P}{\rho } &=&C\gamma \rho ^{\gamma -2}\nabla \rho =:q\left(
\rho \right) \nabla \rho .  \notag
\end{eqnarray}%
Or equivalently%
\begin{equation}
p_{t}+v\cdot \nabla p_{0}+\gamma p_{0}\nabla \cdot v=0  \label{p}
\end{equation}%
We now consider perturbations $\left( \sigma ,v,B,p\right) $ around such a
steady state $\left( \rho _{0},B_{0},V_{0}\equiv 0,P_{0}\right) $ of the
form:%
\begin{align*}
v_{1}& =v_{1}\left( t,x,y\right) ,~v_{2}=v_{2}\left( t,x,y\right) ,\
v_{3}=v_{3}\left( t,x,y\right) , \\
B_{1}& =B_{1}\left( t,x,y\right) ,\ B_{2}=B_{2}\left( t,x,y\right) ,\
B_{3}=B_{3}\left( t,x,y\right) .
\end{align*}%
Equations for perturbed quantities take the form:%
\begin{align}
\sigma _{t}+\nabla \cdot \left[ \left( \rho _{0}+\sigma \right) v\right] &
=0,  \notag \\
\left( \rho _{0}+\sigma \right) \left( v_{t}+v\cdot \nabla v\right) &
=\left( \nabla \times B_{0}\right) \times B+\left( \nabla \times B\right)
\times B_{0}  \label{full-eqn} \\
& +\left( \nabla \times B\right) \times B-\nabla p+\sigma \vec{g},  \notag \\
B_{t}& =\nabla \times \left( v\times \left( B_{0}+B\right) \right) ,  \notag
\\
\nabla \cdot B& =0.  \notag
\end{align}%
We also obtain its linearized system:%
\begin{align}
\sigma _{t}+\nabla \cdot \left( \rho _{0}v\right) & =0,
\label{linearized eqn} \\
\rho _{0}v_{t}& =\left( \nabla \times B_{0}\right) \times B+\left( \nabla
\times B\right) \times B_{0}-\nabla p+\sigma \vec{g},  \notag \\
B_{t}& =\nabla \times \left( v\times B_{0}\right) ,  \notag \\
\nabla \cdot B& =0.  \notag
\end{align}%
In the incompressible case, $v$ also satisfies, in both linear and nonlinear
system,%
\begin{equation*}
\nabla \cdot v=0.\ 
\end{equation*}%
It is crucial to make an equivalent second-order linearized system in one
quantity $v$ and use it through our variational formulation: By taking $t$%
-derivative of $v$-equation in (\ref{linearized eqn}) and plugging $\sigma $%
- and $B$-equation in (\ref{linearized eqn}) (also (\ref{p}) in the
compressible case) into the resulting equation, we obtain the following
second-order linear hyperbolic PDE for velocity $v$:

\begin{equation}
\rho _{0}v_{tt}=\left( Q\cdot \nabla \right) B_{0}+\left( B_{0}\cdot \nabla
\right) Q-\nabla \left( B_{0}\cdot Q+p_{t}\right) -\vec{g}\nabla \cdot
\left( \rho _{0}v\right) =: L\left( v\right) ,  \label{v_tt}
\end{equation}%
where $Q=\nabla \times \left( v\times B_{0}\right) $ and we have used (\ref%
{vec1}).

For notational convenience,

\begin{notation}
For any $u$ and $v\in L^{2}\left( 0,2\pi \right) ,$%
\begin{equation*}
<u,v>=\int_{0}^{2\pi }u\cdot v~dx,\ \left( u,v\right) =\int_{0}^{2\pi }\rho
_{0}u\cdot v~dx,
\end{equation*}%
For any $u$ and $v\in L^{2}\left( D\right) ,$%
\begin{equation*}
<u,v>=\iint_{D}u\cdot v~dxdy,\ \left( u,v\right) =\iint_{D}\rho _{0}u\cdot
v~dxdy.
\end{equation*}
\end{notation}

\begin{notation}
$\left\Vert f\right\Vert =\left( f,f\right) ^{1/2},\ \left\Vert f\right\Vert
_{s}=\left(\sum_{\alpha}{\left\Vert \partial ^{\alpha }f\right\Vert^{2} }%
\right)^{1/2},$ where $\alpha $ is a multi-index with $\left\vert \alpha
\right\vert \leq s.$
\end{notation}

We make the following variational formulation and obtain a discrete set of
eigenvalues for the linear operator (\ref{v_tt}): for each wave number $k\in 
{\mathbb {N}}%
$,%
\begin{equation}
\lambda _{k}^{2}=\sup_{\substack{ u=\left( u_{1},u_{2}\right)  \\ u_{1}\in
H^{\kappa }\left( 0,2\pi \right) ,  \\ u_{2}\in H^{\kappa -1}\left( 0,2\pi
\right) }}\frac{<L\left( u\right) ,u~>}{\left( u,u\right) },  \label{l_k}
\end{equation}%
where $\partial _{y}$ is replaced with multiplying by $k$ in $L\left(
u\right) $ and $\kappa =1$ or $2$. Notice that the RHS of (\ref{l_k}) is
indeed a function of $k$.

In the incompressible case, $u_{2}$ is replaced by $-\frac{1}{k}u_{1x}$ from
divergence-free condition for the normal modes and (\ref{l_k}) reduces to a
variational problem for $u_{1}$ alone. On the other hand, the compressible
case may not be simplified to a formula for $u_{1}$ alone.

Key step of this article is to show that this discrete set $\left\{ \lambda
_{k}\right\} _{k\in 
{\mathbb {N}}%
}$ of eigenvalues characterizes the continuum spectral radius of the whole
linearized operator by taking the limit as $k\rightarrow \infty $ to obtain
the bounded least upper bound $\Lambda >0$: Let%
\begin{equation}
\lim_{k}\lambda _{k}^{2}=\Lambda ^{2}=\sup_{v\in H^{\kappa -1}\left(
D\right) }\frac{\iint_{D}H\left( v,v_{x}\right) dxdy}{\left( v,v\right) },
\label{radius}
\end{equation}%
where coefficients of $H$ depend on $\rho _{0},p_{0},B_{0}.$ Then we obtain
the following theorem:

\begin{theorem}
Let $\left( \sigma ,v\right) $ be a solution to (\ref{linearized eqn}) and
let $\Lambda ^{2}>0$, then we have%
\begin{equation*}
\left\Vert \sigma \left( t\right) ,v\left( t\right) ,B\left( t\right)
\right\Vert _{s}\leq Ce^{\Lambda t}\left\Vert \sigma \left( 0\right)
,v\left( 0\right) ,B\left( 0\right) \right\Vert _{s+2},
\end{equation*}%
where $C=C\left( \Lambda ,\rho _{0},s\right) $.
\end{theorem}

Thanks to this theorem, we can locate a dominant eigenvalue and pass to
nonlinear instability.

We establish the following dynamical instability for the fully nonlinear MHD
system around different steady states satisfying $\Lambda ^{2}>0$ as in
Section 4:

\begin{theorem}
Steady states $\left( \rho _{0},\vec{0},B_{0}\right) $ of (\ref{origin}) in (%
\ref{steady-state}) satisfying $\Lambda ^{2}>0$ are indeed nonlinearly
unstable: For any $s$ large, there exists $\varepsilon _{0}>0$, such that
for any small $\delta >0$, there exists a family of classical solutions $%
(\rho ^{\delta }(t,x,y),V^{\delta }(t,x,y),B^{\delta }\left( t,x,y\right) )$
to $\left( \ref{origin}\right) $\ such that%
\begin{equation*}
\left\Vert \rho ^{\delta }(0,\cdot )-\rho _{0}\left( \cdot \right)
\right\Vert _{H^{s}(D)}+\left\Vert V^{\delta }(0,\cdot )\right\Vert
_{H^{s}(D)}+\left\Vert B^{\delta }\left( 0,\cdot \right) -B_{0}\left( \cdot
\right) \right\Vert _{H^{s}\left( D\right) }\leq \delta ,
\end{equation*}%
but for $T^{\delta }=O(\left\vert \ln \delta \right\vert ),$%
\begin{equation*}
\sup_{0\leq t\leq T^{\delta }}\left\{ \left\Vert \rho ^{\delta }(t,\cdot
)-\rho _{0}\left( \cdot \right) \right\Vert _{L^{1}(D)}+\left\Vert V^{\delta
}(t,\cdot )\right\Vert _{L^{1}(D)}+\left\Vert B^{\delta }\left( t,\cdot
\right) -B_{0}\left( \cdot \right) \right\Vert _{L^{1}\left( D\right)
}\right\} \geq \varepsilon _{0}.
\end{equation*}
\end{theorem}

\begin{remark}
This indicates the Kruskal-Schwarzschild instability for incompressible
ideal MHD fluids as $\gamma \rightarrow \infty $ formally.
\end{remark}

\begin{remark}
The instability time $T^{\delta }$ occurs before the possible blow-up time
which is shown in the proof and we measure instability in $L^{1}$.
\end{remark}

Rayleigh-Taylor instability is well known as gravity-driven instability in
fluids when heavy fluid is on top of light one. Linear instability for an
incompressible fluid was first introduced by Rayleigh in 1883 \cite{R}.
Study on linear instability for incompressible ideal MHD system in the
presence of magnetic field lines could be found in \cite{N}, \cite{CH} for
instance, which includes the classical Kruskal-Schwarzschild instability in
the presence of a magnetic field orthogonal to the gravitational force. In
this case, instability criterion (\ref{rho_0}) and the growth rates (\ref%
{lamkk}), (\ref{Lambda}) turn out to be the same as the ones in the
Rayleigh-Taylor instability as in \cite{HG} without effect of magnetic field
lines on the instability. On the other hand, when we consider effect of a
magnetic field parallel to direction of the gravitational force, stabilizing
effect of magnetic field lines appears as in (\ref{lambdak}) and (\ref%
{Lambda1}). In the presence of a vertical magnetic field, the rigidity
produced by the magnetic field lines hinders its own way to instability and
makes the growth rate slower. Condition for linear instability for a
compressible fluid in the absence of a magnetic field was first derived by
Schwarzschild in \cite{S}, and since been discussed by many other physicists
for a certain class of steady states \cite{CH},\cite{CH2},\cite{Cow},\cite{J}%
. The full consideration of gravity, magnetic field lines and
compressibility has been also largely discussed for its linear instability
by many physicists such as in \cite{B},\cite{BFKK},\cite{P}, which exhibit
very interesting phenomena although it accelerates analytical difficulties.

Despite extensive research and interest in this subject from physical point
of view, little has been done from mathematical perspective for the MHD
system. In addition, the passage from linear to nonlinear instability in a
conservative PDE system is quite difficult because of the following two main
obstacles: (1) presence of the continuum linear spectrum and (2) severely
unbounded high-order terms in PDE systems. No systematic framework has been
built up for this problem although there have been works towards this
subject for specific physical systems, for example \cite{GS},\cite{CGG},\cite%
{HG}. Variational approach was first introduced by Guo and Hwang \cite{HG}
in the case of dynamical Rayleigh-Taylor instability for incompressible
Euler fluids. However it is not obvious whether it can be extended to MHD
instabilities for compressible as well as incompressible fluids since MHD
has more complicated structure in addition to the analytical difficulties
coming from compressibility.

Crucial point is whether and how to locate a dominant eigenvalue in the
complex linear spectrum of the MHD system. We use the MHD energy principle
in order to estimate a sharp spectral radius and make extensive use of the
variational structure of the linearized MHD system, resulting in more
precise and optimal estimates. The radius (\ref{radius}) of continuum
spectrum is obtained as the least upper bound for a discrete set of
eigenvalues (\ref{l_k}) of normal growing modes and this method is explicit
and constructive. We consider three different magnetic fields as steady
states which give rise to different outcome in the growth rates (\ref{Lambda}%
), (\ref{Lambda2}), (\ref{Lambda1}) of instabilities and different ranges of
admissible density profiles for instabilities. We justify stabilizing and
destabilizing effects of magnetic filed lines and compressibility as
expected physically. Furthermore incompressible case can also be viewed as
the limiting case of compressible case as $\gamma \rightarrow \infty $. The
article is organized as follows.

We formulate the variational problems (\ref{l_k}) in Section 2 and show the
existence of smooth maximizers satisfying the corresponding Sturm-Liouville
equations. We then derive a sharp growth rate $\Lambda >0$ for the whole
linear system in Section 3 as in the Theorem 1. In Section 4, we give
different examples of steady states which result in different growth rates,
Sturm-Liouville equations, and different admissibility for density profiles
for instability. For instance, if a steady magnetic field is parallel to the
gravity, we have the fourth-order Sturm-Liouville equation (\ref{SL1}). The
compressible case (\ref{lambda-k}) is more complicated and its variational
problem should be treated with more delicacy.

In Section 5, we construct approximate solutions and give energy estimates
for the full system. Finally, we present nonlinear instability for the ideal
MHD system around our three kind steady states in Section 6.

We put some useful vector identities which are used throughout the paper.
For any two vectors $\vec{a}$ and $\vec{b}$,

\begin{equation}
\nabla \times \left( \vec{a}\times \vec{b}\right) =\vec{a}\left( \nabla
\cdot \vec{b}\right) +\left( \vec{b}\cdot \nabla \right) \vec{a}-\vec{b}%
\left( \nabla \cdot \vec{a}\right) -\left( \vec{a}\cdot \nabla \right) \vec{b%
},  \label{vec1}
\end{equation}%
\begin{equation}
\nabla \left( \vec{a}\cdot \vec{b}\right) =\left( \vec{a}\cdot \nabla
\right) \vec{b}+\left( \vec{b}\cdot \nabla \right) \vec{a}+\vec{a}\times
\left( \nabla \times \vec{b}\right) +\vec{b}\times \left( \nabla \times \vec{%
a}\right) ,  \label{vec2}
\end{equation}

\section{General variational framework}

We consider the following steady magnetic fields:

(Case $B_{0}\perp \vec{g}$)

\begin{equation*}
B_{0}=\left( 0,0,B_{0}\left( x\right) \right) ,
\end{equation*}

(Case $B_{0}\parallel \vec{g}$)%
\begin{equation*}
B_{0}=\left( B_{0},0,0\right) .
\end{equation*}%
In the case of $B_{0}\parallel \vec{g}$, we assume $v_{3}=B_{3}=0$. By
integration by parts, using (\ref{linearized eqn}),(\ref{p}),(\ref%
{steady-state}) and completing the square with respect to $\nabla \cdot v$,
we obtain the following decomposition for $<L\left( v\right) ,v>$: for any $%
v=\left( v_{1},v_{2},v_{3}\right),$%
\begin{eqnarray}
&<&L\left( v\right) ,v>=\iint_{D}L\left( v\right) \cdot v\ dxdy  \label{L} \\
&=&\iint_{D}\left[ -F\left( v_{2x}\right) -G\left( \nabla \cdot
v,v_{1},v_{1x}\right) +H\left( v_{1},v_{1x}\right) \right] dxdy.  \notag
\end{eqnarray}%
We state some important properties which are satisfied by the above
functional:

\begin{enumerate}
\item $L$ is variational, i.e., for any $u$ and $v$%
\begin{equation*}
<L\left( u\right) ,v>=<u,L\left( v\right) >.
\end{equation*}

\item $<L\left( v\right) ,v>$ is concave with respect to $%
v_{1x},v_{2x},v_{2} $ respectively.

\item $<L\left( v\right) ,v>$ $\rightarrow -\infty $ as $v_{1x}\rightarrow
\infty $ and $<L\left( v\right) ,v>$ $\rightarrow -\infty $ as $%
v_{2x}\rightarrow \infty $ if $F\neq 0$.

\item $F>0,G>0,$ $H\left( 0\right) =0$ when $H=H\left( v_{1}\right) $ alone,
and $F,G,H$ are all quadratic.
\end{enumerate}

\begin{remark}
$G=0$ in the incompressible case.
\end{remark}

\begin{remark}
For the case $B_{0}\perp \vec{g}$, we have $F=0$, $H=H\left( v_{1}\right) ,$
and $G=G\left( \nabla \cdot v,v_{1}\right) $ with $b=0$ while we have $%
F=F\left( v_{2x}\right) ,H=H\left( v_{1},v_{1x}\right) ,G=G\left( \nabla
\cdot v,v_{1},v_{1x}\right) $ in the case $B_{0}\parallel \vec{g}$ as we can
see in Section 4.
\end{remark}

We will show that for any fixed wave number $k\in 
{\mathbb {N}}%
$, the corresponding eigenvalue $\lambda _{k}>0$ for the linearized MHD
system takes the variational formulation closely related to the above
variational structure. A normal mode is of the form:%
\begin{eqnarray}
v_{1}\left( t,x,y\right) &=&\tilde{v}_{1}\left( x\right) \cos \left(
ky\right) \exp \left( \lambda _{k}t\right) ,\   \label{expo} \\
v_{2}\left( t,x,y\right) &=&\tilde{v}_{2}\left( x\right) \sin \left(
ky\right) \exp \left( \lambda _{k}t\right) ,  \notag \\
v_{3}\left( t,x,y\right) &=&\tilde{v}_{3}\left( x\right) \cos \left(
ky\right) \exp \left( \lambda _{k}t\right) ,  \notag
\end{eqnarray}%
where $k$ is a wave number. Substituting (\ref{expo}) into (\ref{v_tt})
yields the following second-order ODE for $\tilde{v}$:%
\begin{equation}
\lambda _{k}^{2}\rho _{0}\left( \tilde{v}_{1},\tilde{v}_{2},\tilde{v}%
_{3}\right) =L\left( \tilde{v}_{1},\tilde{v}_{2},\tilde{v}_{3}\right) .
\label{SL}
\end{equation}%
Note that $L_{3}\left( \tilde{v}_{1},\tilde{v}_{2},\tilde{v}_{3}\right) =0$
and thus $\tilde{v}_{3}=0$. We now make the following variational
formulations:%
\begin{equation}
\lambda _{k}^{2}=\sup_{\substack{ u=\left( u_{1},u_{2}\right)  \\ 
_{\substack{ u_{1}\in H^{\kappa }\left( 0,2\pi \right)  \\ u_{2}\in
H^{\kappa -1}\left( 0,2\pi \right) }}}}\frac{<L\left( u\right) ,u~>}{\left(
u,u\right) },  \label{lam}
\end{equation}%
where $\kappa =1$ or $2$, $\partial _{y}$ is replaced with multiplying by $k$
in $L\left( u\right) $ and integrations here are with respect to $x$ over $%
\left( 0,2\pi \right) $. Indeed,

Incompressible case: we use divergence-free condition to reduce to%
\begin{equation*}
\lambda _{k}^{2}=\sup_{u\in H^{\kappa }\left( 0,2\pi \right) }\frac{%
\int_{0}^{2\pi }\left[ -F\left( -\frac{u_{xx}}{k}\right) +H\left(
u,u_{x}\right) \right] dx}{\int_{0}^{2\pi }\left[ u^{2}+\frac{u_{x}^{2}}{%
k^{2}}\right] dx},
\end{equation*}%
where $\kappa =1$ or $2$.

Compressible case:%
\begin{equation*}
\lambda _{k}^{2}=\sup_{\substack{ u_{1}\in H^{\kappa }\left( 0,2\pi \right) 
\\ u_{2}\in H^{\kappa -1}\left( 0,2\pi \right) }}\frac{\int_{0}^{2\pi }\left[
-F\left( u_{2x}\right) -G\left( u_{1x}+ku_{2},u_{1},u_{1x}\right) +H\left(
u_{1},u_{1x}\right) \right] dx}{\int_{0}^{2\pi }\left[ u_{1}^{2}+u_{2}^{2}%
\right] dx},
\end{equation*}%
where $\kappa =1$ or $2$.

Let 
\begin{equation}
\Lambda ^{2}=\sup_{v\in H^{\kappa -1}\left( D\right) }\frac{\iint_{D}H\left(
v,v_{x}\right) dxdy}{\left( v,v\right) },  \label{Lam}
\end{equation}%
where $\kappa =1$ or $2$. Then we show $\Lambda^2 $ is the least upper bound
for $\left\{ \lambda _{k}^2\right\} _{k\in N}$.

\begin{lemma}
\begin{equation*}
\lim_{k\rightarrow \infty }\lambda _{k}^{2}=\Lambda ^{2}.
\end{equation*}
\end{lemma}

\begin{proof}
Since $F\left( 0\right) =0$, it is easy to see that our Lemma is true for
incompressible case by letting $k\rightarrow \infty $. We now treat
compressible case.

Note that, with the choice of $u_{2}=-\frac{1}{k}\left(
au_{1}+bu_{1x}\right) $, we have%
\begin{eqnarray*}
\lambda _{k}^{2} &=&\sup_{\substack{ u=\left( u_{1},u_{2}\right)  \\ 
_{\substack{ u_{1}\in H^{\kappa }\left( 0,2\pi \right)  \\ u_{2}\in
H^{\kappa -1}\left( 0,2\pi \right) }}}}\frac{<L\left( u\right) ,u>}{\left(
u,u\right) } \\
&\geq &\sup_{u_{1}\in H^{\kappa }\left( 0,2\pi \right) }\frac{\int_{0}^{2\pi
}\left[ -F\left( -\frac{1}{k}\left( au_{1}+bu_{1x}\right) _{x}\right)
+H\left( u_{1},u_{1x}\right) \right] dx}{\int_{0}^{2\pi }\rho _{0}\left[
u_{1}^{2}+\frac{1}{k^{2}}\left( au_{1}+bu_{1x}\right) ^{2}\right] dx}.
\end{eqnarray*}%
Thus, we have the following inequality%
\begin{equation*}
\sup_{u_{1}\in H^{\kappa }\left( 0,2\pi \right) }\frac{\int_{0}^{2\pi }\left[
-F\left( -\frac{1}{k}\left( au_{1}+bu_{1x}\right) _{x}\right) +H\left(
u_{1},u_{1x}\right) \right] dx}{\int_{0}^{2\pi }\rho _{0}\left[ u_{1}^{2}+%
\frac{1}{k^{2}}\left( au_{1}+bu_{1x}\right) ^{2}\right] dx}\leq \lambda
_{k}^{2}\leq
\end{equation*}%
\begin{eqnarray*}
&&\sup_{\substack{ u_{1}\in H^{\kappa }\left( 0,2\pi \right)  \\ u_{2}\in
H^{\kappa -1}\left( 0,2\pi \right) }}\frac{\int_{0}^{2\pi }\left[ -F\left(
u_{2x}\right) -G\left( u_{1x}+ku_{2},u_{1},u_{1x}\right) \right] dx}{%
\int_{0}^{2\pi }\rho _{0}\left[ u_{1}^{2}+u_{2}^{2}\right] dx}+\sup_{ 
_{\substack{ u_{1}\in H^{\kappa }\left( 0,2\pi \right)  \\ u_{2}\in
H^{\kappa -1}\left( 0,2\pi \right) }}}\frac{\int_{0}^{2\pi }H\left(
u_{1},u_{1x}\right) dx}{\int_{0}^{2\pi }\rho _{0}\left[ u_{1}^{2}+u_{2}^{2}%
\right] dx} \\
&\leq &\sup_{u_{1}\in H^{\kappa -1}\left( 0,2\pi \right) }\frac{%
\int_{0}^{2\pi }H\left( u_{1},u_{1x}\right) dx}{\int_{0}^{2\pi }\rho
_{0}u_{1}^{2}dx}.
\end{eqnarray*}%
Thanks to%
\begin{equation*}
\frac{1}{k}\left( au_{1}+bu_{1x}\right) \rightarrow 0~\text{as }k\rightarrow
\infty ,\ F\left( 0\right) =0,
\end{equation*}%
letting $k\rightarrow \infty $ yields%
\begin{equation*}
\lim_{k\rightarrow \infty }\lambda _{k}^{2}=\sup_{u_{1}\in H^{\kappa
-1}\left( 0,2\pi \right) }\frac{\int_{0}^{2\pi }H\left( u_{1},u_{1x}\right)
dx}{\int_{0}^{2\pi }\rho _{0}u_{1}^{2}dx},
\end{equation*}%
and%
\begin{equation*}
\lim_{k\rightarrow \infty }\lambda _{k}^{2}=\sup_{u_{1}\in H^{\kappa
-1}\left( D\right) }\frac{\iint_{D}H\left( u_{1},u_{1x}\right) dxdy}{%
\iint_{D}\rho _{0}u_{1}^{2}dxdy},
\end{equation*}%
since $H^{\kappa }$ is dense in $H^{\kappa -1}$. Thus the proof is complete.
\end{proof}

We now show the existence of maximizer for the variational problem (\ref{lam}%
). Assume that%
\begin{equation}
\left( u,u\right) =\int_{0}^{2\pi }\rho _{0}\left[ u_{1}^{2}+u_{2}^{2}\right]
~dx=1.  \label{constraint-k}
\end{equation}%
For fixed $k,$ let

\begin{equation}
\lambda _{k}^{2}=\sup_{_{\substack{ u_{1}\in H^{\kappa }\left( 0,2\pi
\right)  \\ u_{2}\in H^{\kappa -1}\left( 0,2\pi \right) }}}<L\left( u\right)
,u>.  \label{va}
\end{equation}

\begin{lemma}
For any fixed $k,$ there exists a smooth maximizer for the variational
problem (\ref{va}) with the constraint (\ref{constraint-k}).
\end{lemma}

\begin{proof}
Let $\left\{ u_{1}^{n},u_{2}^{n}\right\} $ be a maximizing sequence with the
constraint (\ref{constraint-k}). Then $u_{1}^{n}$ and $u_{2}^{n}$ converge
weakly in $L^{2}\left( 0,2\pi \right) $ to $u_{1}^{0}$ and $u_{2}^{0}$
respectively and we have%
\begin{equation}
<L\left( u^{n}\right) ,u^{n}>\rightarrow \lambda _{k}^{2}.  \label{conv}
\end{equation}

\textbf{Case 1} $B_{0}\perp \vec{g}$, where $F=0,H=H\left( u_{1}\right) $:

Since $<L\left( u\right) ,u>\rightarrow -\infty $ as $u_{1x}\rightarrow
\infty $ (Property 3 of the functional), $u_{1x}$ is bounded in $L^{2}\left(
0,2\pi \right) $ uniformly in $n$. Thus there exists a weak limit $\left\{
u_{1}^{0},u_{2}^{0}\right\} $ such that%
\begin{eqnarray*}
u_{1x}^{n}\text{ } &\rightharpoonup &\text{ }u_{1x}^{0}\text{ weakly in }%
L^{2}\left( 0,2\pi \right) ,\text{ }u_{2x}^{n}\rightharpoonup \text{ }%
u_{2}^{0}\text{ weakly in }L^{2}\left( 0,2\pi \right) \\
u_{1}^{n} &\rightarrow &u_{1}^{0}\text{ strongly in }L^{2}\left( 0,2\pi
\right) \text{ and }u_{1}^{0}\in H^{1}\left( 0,2\pi \right) .
\end{eqnarray*}%
Next we show that $\left\{ u_{1}^{0},u_{2}^{0}\right\} $ is a maximizer and
satisfies the constraint (\ref{constraint-k}). Since $\left\{
u_{1}^{0},u_{2}^{0}\right\} $ is a weak limit of $\left\{
u_{1}^{n},u_{2}^{n}\right\} $, $\left( u^{0},u^{0}\right) =\int_{0}^{2\pi
}\rho _{0}\left[ \left( u_{1}^{0}\right) ^{2}+\left( u_{2}^{0}\right) ^{2}%
\right] dx\leq 1$ by lower semi-continuity of $L^{2}$:%
\begin{eqnarray*}
&&\int_{0}^{2\pi }\rho _{0}\left[ \left( u_{1}^{n}\right) ^{2}+\left(
u_{2}^{n}\right) ^{2}\right] dx-\int_{0}^{2\pi }\rho _{0}\left[ \left(
u_{1}^{0}\right) ^{2}+\left( u_{2}^{0}\right) ^{2}\right] dx \\
&=&\int_{0}^{2\pi }2\rho _{0}\left[ u_{1}^{0}\left(
u_{1}^{n}-u_{1}^{0}\right) +u_{2}^{0}\left( u_{2}^{n}-u_{2}^{0}\right) %
\right] dx \\
&&+\int_{0}^{2\pi }\rho _{0}\left[ \left( u_{1}^{n}-u_{1}^{0}\right)
^{2}+\left( u_{2}^{n}-u_{2}^{0}\right) ^{2}\right] \\
&\geq &\int_{0}^{2\pi }2\rho _{0}\left[ u_{1}^{0}\left(
u_{1}^{n}-u_{1}^{0}\right) +u_{2}^{0}\left( u_{2}^{n}-u_{2}^{0}\right) %
\right] dx\rightarrow 0\text{ as }n\rightarrow 0.
\end{eqnarray*}%
due to the fact that $\left\{ u_{1}^{0},u_{2}^{0}\right\} $ is a weak limit
of $\left\{ u_{1}^{n},u_{2}^{n}\right\} $ in $L^{2}\left( 0,2\pi \right) $.
We use concavity of the functional $<L\left( u\right) ,u>$ with respect to $%
u_{1x},u_{2x},u_{2}$ (Property 2 of the functional) and the strong
convergence of $u_{1}^{n}$ to $u_{1}^{0}$ in $L^{2}\left( 0,2\pi \right) $
to deduce%
\begin{equation}
<L\left( u^{0}\right) ,u^{0}>\ \geq \lambda _{k}^{2}.  \label{conc}
\end{equation}%
Let 
\begin{equation*}
J\left( u_{2x},u_{1x},u_{2},u_{1}\right) =<L\left( u\right) ,u>.
\end{equation*}%
Then by concavity of this functional \thinspace $J$, strong convergence of $%
u_{1}^{n}$ to $u_{1}^{0}$ and weak convergence of $%
u_{1x}^{n},u_{2}^{n},u_{1}^{n}$ to $u_{1x}^{0},u_{2}^{0},u_{1}^{0}$
respectively, we have%
\begin{eqnarray*}
&&J\left( u_{2x}^{n},u_{1x}^{n},u_{2}^{n},u_{1}^{n}\right) -J\left(
u_{2x}^{0},u_{1x}^{0},u_{2}^{0},u_{1}^{0}\right) \\
&=&J\left( u_{2x}^{n},u_{1x}^{n},u_{2}^{n},u_{1}^{n}\right) -J\left(
u_{2x}^{n},u_{1x}^{n},u_{2}^{n},u_{1}^{0}\right) +J\left(
u_{2x}^{n},u_{1x}^{n},u_{2}^{n},u_{1}^{0}\right) -J\left(
u_{2x}^{0},u_{1x}^{0},u_{2}^{0},u_{1}^{0}\right) \\
&\leq &J\left( u_{2x}^{n},u_{1x}^{n},u_{2}^{n},u_{1}^{n}\right) -J\left(
u_{2x}^{n},u_{1x}^{n},u_{2}^{n},u_{1}^{0}\right) +\nabla J\left(
u_{2x}^{0},u_{1x}^{0},u_{2}^{0},u_{1}^{0}\right) \cdot \left(
u_{2x}^{n}-u_{2x}^{0},u_{1x}^{n}-u_{1x}^{0},u_{2}^{n}-u_{2}^{0},0\right) \\
&\rightarrow &0\text{ \ as }n\rightarrow \infty .
\end{eqnarray*}%
Thus we obtain (\ref{conc}). Notice that we can not have $u^{0}=\left\{
u_{1}^{0},u_{2}^{0}\right\} =\left\{ 0,0\right\} $ a.e. since, by strong
convergence of $u_{1}^{n}$ to $u_{1}^{0}$ and by (\ref{conv}), we have 
\begin{equation*}
\int_{0}^{2\pi }H\left( u_{1}^{0}\right) dx=\lim_{n\rightarrow \infty
}\int_{0}^{2\pi }H\left( u_{1}^{n}\right) dx\geq \lambda _{k}^{2}>0,
\end{equation*}%
where we have used the property 4 of the functional. Suppose now that $%
\left( u^{0},u^{0}\right) =\alpha ^{2}<1.$ Then let $\left( \tilde{u}%
_{1}^{0},\tilde{u}_{2}^{0}\right) =\frac{1}{\alpha }\left(
u_{1}^{0},u_{2}^{0}\right) $ so that $\left( \tilde{u}^{0},\tilde{u}%
^{0}\right) =1$. By the above argument (\ref{conc}), we have%
\begin{equation*}
<L\left( \tilde{u}^{0}\right) ,\tilde{u}^{0}>\ \geq \frac{\lambda _{k}^{2}}{%
\alpha }>\lambda _{k}^{2},
\end{equation*}%
leading to a contradiction. Thus $\left\{ u_{1}^{0},u_{2}^{0}\right\} $ is a
maximizer satisfying the constraint (\ref{constraint-k}).

\textbf{Case 2 }$B_{0}\parallel \vec{g}$, where $F\left( u_{2x}\right)
=B_{0}^{2}u_{2x}^{2}$:

In this case, by property 3 of the functional, $u_{2x}^{n}$ is also bounded
in $L^{2}\left( 0,2\pi \right) $ uniformly in $n$ and thus both $u_{1}^{n}$
and $u_{2}^{n}$ converge strongly in $L^{2}\left( 0,2\pi \right) $ to $%
u_{1}^{0}$ and $u_{2}^{0}$ respectively. Then we have%
\begin{equation*}
\left( u^{0},u^{0}\right) =\int_{0}^{2\pi }\rho _{0}\left[ \left(
u_{1}^{0}\right) ^{2}+\left( u_{2}^{0}\right) ^{2}\right] dx=1.
\end{equation*}%
In a similar manner, we obtain (\ref{conc}) and hence $u^{0}$ is a maximizer.

We finally show such a maximizer satisfies the generalized Sturm-Liouville
equation (\ref{SL}) for both cases. For $\tau \in 
{\mathbb R}%
$ and $w=\left\{ w_{1},w_{2}\right\} \in H^{\kappa }\left( 0,2\pi \right)
\times H^{\kappa -1}\left( 0,2\pi \right) ,$ define $u\left( \tau \right)
=u^{0}+\tau w$, then by (\ref{lambda-k}), we have%
\begin{equation*}
<L\left( u\left( \tau \right) \right) ,u\left( \tau \right) >\ \leq \lambda
_{k}^{2}\ \left( u\left( \tau \right) ,u\left( \tau \right) \right) .
\end{equation*}%
Set%
\begin{equation*}
I\left( \tau \right) =<L\left( u\left( \tau \right) \right) ,u\left( \tau
\right) >-\lambda _{k}^{2}\ \left( u\left( \tau \right) ,u\left( \tau
\right) \right) ,
\end{equation*}%
then we have $I\left( \tau \right) \leq 0$ for all $\tau \in 
{\mathbb R}%
$ and $I\left( 0\right) =0$. This implies%
\begin{eqnarray*}
I^{\prime }\left( 0\right) &=&<L\left( w\right) ,u^{0}>+<L\left(
u^{0}\right) ,w>-2\lambda _{k}^{2}\left( u^{0},w\right) \\
&=&2<L\left( u^{0}\right) -\lambda _{k}^{2}\rho _{0}u^{0},w>=0\text{ for all 
}w,
\end{eqnarray*}%
since $L$ is variational (Property 1 of the functional). Thus $u^{0}$
satisfies a normal mode, i.e.,%
\begin{equation*}
\lambda _{k}^{2}\rho _{0}u^{0}=L\left( u^{0}\right) .
\end{equation*}%
Since $\rho _{0}$, $p_{0}$, and $B_{0}$ are smooth, $u^{0}$ is also smooth.
This completes the proof.
\end{proof}

\section{Linear growth rate $\Lambda $}

In this section, we show that $\Lambda $ is the optimal growth rate for the
linearized system and it serves as the spectral radius of the linear
operator. We state a global existence of solutions to the linearized system,
which can be obtained by a straightforward method.

\begin{lemma}
There exists a global in time solution $\left( \sigma ,v,B\right) \in
C\left( \left[ 0,T\right] ;H^{s}\left( D\right) \right) $ to the linearized
MHD system (\ref{linearized eqn}).
\end{lemma}

Let $\left( \sigma ,v\right) $ be a solution to (\ref{linearized eqn}) and
let $\Lambda ^{2}>0$, then we have

\begin{theorem}
\begin{equation*}
\left\Vert \sigma \left( t\right) ,v\left( t\right) ,B\left( t\right)
\right\Vert _{s}\leq Ce^{\Lambda t}\left\Vert \sigma \left( 0\right)
,v\left( 0\right) ,B\left( 0\right) \right\Vert _{s+2},
\end{equation*}%
where $C=C\left( \Lambda ,\rho _{0},s\right) $ and $\Lambda >0$.
\end{theorem}

\begin{proof}
We show by induction on $\bar{s}$, the number of $x$-derivatives. We first
treat the case $\bar{s}=0$: Multiply (\ref{2nd-linear eq}) by $v_{t}$ and
integrate over $x$ and $y$, then we have%
\begin{equation}
\frac{d}{dt}\left( v_{t},v_{t}\right) =\frac{d}{dt}<L\left( v\right) ,v>
\label{v-t}
\end{equation}%
since $L$ is variational. Notice that using (\ref{Lam}) yields%
\begin{equation*}
\iint_{D}H\left( v_{1},v_{1x}\right) dxdy\leq \Lambda ^{2}\left(
v_{1},v_{1}\right) \leq \Lambda ^{2}\left( v,v\right) .
\end{equation*}%
By integrating (\ref{v-t}) over time and by (\ref{L}), we obtain%
\begin{eqnarray}
&&\left( v_{t},v_{t}\right) +\iint_{D}\left[ F\left( v_{2x}\right) +G\left(
\nabla \cdot v,v_{1},v_{1x}\right) \right] dxdy  \notag \\
&\leq &\iint_{D}H\left( v_{1},v_{1x}\right) dxdy+\left\Vert \sigma \left(
0\right), v\left( 0\right), B\left( 0\right) \right\Vert _{1}^{2}.  \notag \\
&\leq &\Lambda ^{2}\left( v,v\right) +\left\Vert \sigma \left( 0\right),
v\left( 0\right), B\left( 0\right) \right\Vert _{1}^{2},  \label{div}
\end{eqnarray}%
and%
\begin{equation*}
\frac{d}{dt}\left\Vert v\right\Vert \leq \left\Vert v_{t}\right\Vert \leq
\left\Vert \sigma \left( 0\right) ,v\left( 0\right), B\left( 0\right)
\right\Vert _{1}+\Lambda \left\Vert v\right\Vert .
\end{equation*}%
Thus, we have%
\begin{equation*}
\left\Vert v\right\Vert \leq Ce^{\Lambda t}.
\end{equation*}%
where $C=\left\Vert \sigma \left( 0\right) ,v\left( 0\right), B\left(
0\right) \right\Vert _{1} $ and from now on, we will just use universal
constant $C$ which varies, only for notational convenience. Notice that
these estimates exactly apply to the $t$- and $y$-derivatives of any order
of $v$ and $\sigma $ since the variational structure of (\ref{v_tt}) is not
destroyed by taking $t$- and $y$-derivatives.

\begin{case}
$B_{0}\perp \vec{g}$ where $F=0,\ G=\left( \gamma p_{0}+B_{0}^{2}\right)
\left( \nabla \cdot v+\frac{g\rho _{0}}{\gamma p_{0}+B_{0}^{2}}v_{1}\right)
^{2}$ ($G=0$ in the incompressible case)$,\ H=H\left( v_{1}\right) ,\ \left(
Q\cdot \nabla \right) B_{0}+\left( B_{0}\cdot \nabla \right) Q=0$:
\end{case}

By (\ref{div}) and (\ref{linearized eqn}), we also have%
\begin{equation*}
\left\Vert \nabla \cdot v\right\Vert \leq Ce^{\Lambda t},\left\Vert \sigma
\right\Vert \leq Ce^{\Lambda t},\left\Vert B\right\Vert \leq Ce^{\Lambda t}.
\end{equation*}%
Next, we consider the case $\bar{s}=1$. By taking curl $v$-equation of (\ref%
{linearized eqn}), we have, with $\omega =\nabla \times v$,%
\begin{equation}
\rho _{0}\left( \omega \right) _{tt}=-\nabla \rho _{0}\times v_{tt}-g\sigma
_{ty}\hat{k},  \label{curl}
\end{equation}%
where $\hat{k}$ is the unit vector in the $z$-direction. Since $\sigma _{ty}$
and $v_{tt}$ have no $x$-derivatives and $B_{1x}=-B_{2y}$ and so have the
growth rate $\Lambda $ as in the previous step, we have%
\begin{equation*}
\left\Vert \omega \right\Vert \leq Ce^{\Lambda t},
\end{equation*}%
where $C=\left\Vert \sigma \left( 0\right) ,v\left( 0\right) ,B\left(
0\right) \right\Vert _{3}$. Thanks to the identity $\Delta \zeta =-\nabla
\times \left( \nabla \times \zeta \right) +\nabla \left( \nabla \cdot \zeta
\right) $ for any $\zeta =\left( \zeta _{1},\zeta _{2}\right) $, we conclude
that all the first derivatives of $v$ have the same growth rate $\Lambda $.
Now for $\nabla \sigma $, we use the vector identity (\ref{vec1}) to get 
\begin{equation}
B_{t}=\nabla \left( v\times B_{0}\right) =-B_{0x}v_{1}\hat{k}-\hat{k}%
B_{0}\nabla \cdot v.  \label{b}
\end{equation}%
Plugging (\ref{b}) and (\ref{p}) into (\ref{v_tt}) yields%
\begin{equation}
\rho _{0}v_{tt}=\nabla \left( B_{0}B_{0x}v_{1}\right) +\nabla \left( \left[
\gamma p_{0}+B_{0}^{2}\right] \nabla \cdot v\right) +\sigma _{t}\vec{g}.
\label{vtt}
\end{equation}%
By induction hypotheses and $\sigma $- and $B$-equations in (\ref{linearized
eqn}), we deduce 
\begin{equation*}
\left\Vert \nabla \left( \nabla \cdot v\right) \right\Vert \leq Ce^{\Lambda
t},\left\Vert \nabla \sigma \right\Vert \leq Ce^{\Lambda t},\left\Vert
\nabla B\right\Vert \leq Ce^{\Lambda t}.
\end{equation*}%
For higher derivatives when $\bar{s}\geq 2$, we use induction hypotheses.
Suppose all the derivatives of order $<\bar{s}$ have the growth rate $%
\Lambda $. Let $\alpha $ be multi index whose order of $x$-derivative is $%
\bar{s}-1$. Then, by taking $\partial _{\alpha }$ of the curl of (\ref%
{2nd-linear eq}), curl$\left( \partial _{\alpha }v\right) $ satisfies%
\begin{equation}
\rho _{0}\left( \partial _{\alpha }\omega \right) _{tt}=-\partial _{a}\left(
\nabla \times \left( \rho _{0}v\right) \right) _{tt}+\rho _{0}\left(
\partial _{\alpha }\omega \right) _{tt}-g\partial _{\alpha }\sigma _{ty}\hat{%
k}.  \label{linear-induction2}
\end{equation}%
The right hand side of (\ref{linear-induction2}) contains derivatives of $v$
and $\sigma $ whose $x$-order $<\bar{s}$. Clearly the $x$-order of $\partial
_{\alpha }\sigma _{ty}$ is still $\bar{s}-1$. The first and the second term
of the RHS together are of the form 
\begin{equation}
\sum_{_{\substack{ \left\vert \beta \right\vert \geq 1  \\ \beta +\gamma
=\alpha }}}C_{\beta ,\gamma }\nabla \times \left[ \left( \partial _{\beta
}\rho _{0}\right) \left( \partial _{\gamma }v_{tt}\right) \right] .
\label{sum1}
\end{equation}%
Notice that only purely $x$-derivatives $\partial _{\beta }\rho _{0}$ with $%
\left\vert \beta \right\vert \geq 1$ in (\ref{sum1}) remain nonzero since $%
\rho _{0}\left( x\right) $ depends only on $x$. Then (\ref{linear-induction2}%
) yields%
\begin{equation*}
\left\Vert \nabla \times \partial _{a}v\right\Vert \leq Ce^{\Lambda t}.
\end{equation*}%
For $\nabla \cdot \partial _{\alpha }v$, we take $\partial _{\beta }$ of the
first component of (\ref{vtt}) with $\partial _{\beta }\partial
_{x}=\partial _{\alpha }$ to get%
\begin{equation*}
\partial _{\beta }\left( \rho _{0}v_{1tt}\right) =\partial _{\alpha }\left(
B_{0}B_{0x}v_{1}\right) +\partial _{\alpha }\left( \left[ \gamma
p_{0}+B_{0}^{2}\right] \nabla \cdot v\right) +\partial _{\beta }\sigma _{t}g.
\end{equation*}%
Using the induction hypotheses as for the case $\bar{s}=1$, we deduce%
\begin{equation*}
\left\Vert v\right\Vert _{\bar{s}}\leq Ce^{\Lambda t},\left\Vert \sigma
\right\Vert _{\bar{s}}\leq Ce^{\Lambda t},\left\Vert B\right\Vert _{\bar{s}%
}\leq Ce^{\Lambda t},
\end{equation*}%
where $C=\left\Vert \sigma \left( 0\right) ,v\left( 0\right) ,B\left(
0\right) \right\Vert _{\bar{s}+2}$. In the incompressible case, we use (\ref%
{curl}) together with the $\nabla \cdot v=0$ to deduce the Lemma.

\begin{case}
$B_{0}\parallel \vec{g}$ and incompressible case where $%
F=B_{0}^{2}v_{2x}^{2},G=0,H=H\left( v_{1},v_{1x}\right) $:
\end{case}

From (\ref{div}) and div $v=0$ and using the induction hypotheses, we obtain%
\begin{equation*}
\left\Vert v_{2x}\right\Vert \leq Ce^{\Lambda t},\left\Vert
v_{1x}\right\Vert \leq Ce^{\Lambda t}.
\end{equation*}%
In this case, by (\ref{vec1}), we have $Q=B_{0}v_{x}$ and $B_{0}\cdot
Q=B_{0}^{2}v_{1x}$. Then (\ref{v_tt}) yields%
\begin{equation}
\rho _{0}v_{tt}=B_{0}^{2}v_{xx}-\nabla \left( B_{0}\cdot Q+p_{t}\right)
+\sigma _{t}\vec{g}.  \label{v-xx}
\end{equation}%
Taking the curl of (\ref{v-xx}) yields%
\begin{equation*}
\rho _{0}\omega _{tt}=B_{0}^{2}\omega _{xx}-\nabla \rho _{0}\times
v_{tt}-g\sigma _{ty}\hat{k}.
\end{equation*}%
Using the induction hypotheses and div $v=0$, we can deduce the Lemma.
Therefore the proof is complete.
\end{proof}

\section{Examples}

In this section, we present examples of steady states which give rise to
different results.

(1) Magnetic field is transverse to the gravity $B_{0}\perp \vec{g}$ for
incompressible fluids:

\begin{equation}
B_{0}=\left( 0,0,B_{0}\left( x\right) \right) .  \label{B1}
\end{equation}%
Linear operator $L$ reads 
\begin{equation}
L\left( v\right) :=-\nabla \left( p_{t}+B_{0}\cdot Q\right) -\vec{g}\rho
_{0x}v_{1},  \notag
\end{equation}%
Sturm-Liouville equation is%
\begin{equation*}
\rho _{0}\lambda _{k}^{2}\tilde{v}_{1}=\frac{\lambda _{k}^{2}}{k^{2}}\left(
\rho _{0}\tilde{v}_{1x}\right) _{x}-g\rho _{0x}\tilde{v}_{1},
\end{equation*}%
Variational formulation for $\lambda _{k}^{2}$ is 
\begin{equation}
\lambda _{k}^{2}=\sup_{u\in H^{1}\left( 0,2\pi \right) }\frac{\int_{0}^{2\pi
}L\left( u\right) \cdot u\ dx}{\int_{0}^{2\pi }\rho _{0}\left[ u^{2}+\frac{%
u_{x}^{2}}{k^{2}}\right] \,dx}>0,  \label{lamkk}
\end{equation}%
where

\begin{equation*}
L\left( u\right) \cdot u\ =H\left( u\right) =-g\rho _{0x}u^{2}
\end{equation*}%
Spectral radius $\Lambda >0$ is given by

\begin{equation}
\Lambda ^{2}=\sup_{v\in L^{2}\left( D\right) }\frac{\iint_{D}H\left(
v\right) \,dxdy}{\iint_{D}\rho _{0}v^{2}\,dxdy}>0.  \label{Lambda}
\end{equation}%
The sufficient condition for instability in the case (1) is given by \ 
\begin{equation}
\rho _{0x}\left( x_{0}\right) <0\text{ \ \ for some }x_{0}\in \left( 0,2\pi
\right) ,  \label{rho_0}
\end{equation}%
(2) Magnetic field is transverse to the gravity $B_{0}\perp \vec{g}$ for
compressible fluids:%
\begin{equation}
B_{0}=\left( 0,0,B_{0}\left( x\right) \right) .  \label{B3}
\end{equation}%
Linear operator $L$, Sturm-Liouville equation, and $\lambda _{k}$ are as
follows: 
\begin{eqnarray}
L(v) &=&-\nabla \left( B_{0}\cdot Q\right) -\nabla \left( p_{0x}v_{1}+\gamma
p_{0}\nabla \cdot v\right) -\vec{g}\left( \rho _{0x}v_{1}+\rho _{0}\nabla
\cdot v\right)  \label{2nd-linear eq} \\
&&  \notag
\end{eqnarray}
\begin{equation}
\lambda _{k}^{2}\rho _{0}\left( \tilde{v}_{1},\tilde{v}_{2},\tilde{v}%
_{3}\right) =L\left( \tilde{v}_{1},\tilde{v}_{2},\tilde{v}_{3}\right) ,
\label{F-k}
\end{equation}%
where $L$ $=\left( L^{1},L^{2},L^{3}\right) $is a linear functional in $%
\tilde{v}$ $=\left( \tilde{v}_{1},\tilde{v}_{2},\tilde{v}_{3}\right) $ such
that 
\begin{eqnarray*}
L^{1} &=&\left( g\rho _{0}\tilde{v}_{1}+\left( \gamma p_{0}+B_{0}^{2}\right)
\left( \tilde{v}_{1x}+k\tilde{v}_{2}\right) \right) _{x}-g\left( \rho _{0x}%
\tilde{v}_{1}+\rho _{0}\left( \tilde{v}_{1x}+k\tilde{v}_{2}\right) \right) ,
\\
L^{2} &=&-k\left( g\rho _{0}\tilde{v}_{1}+\left( \gamma
p_{0}+B_{0}^{2}\right) \left( \tilde{v}_{1x}+k\tilde{v}_{2}\right) \right) ,
\\
L^{3} &=&0,
\end{eqnarray*}%
and so $\tilde{v}_{3}=0$.

\begin{equation}
\lambda _{k}^{2}=\sup_{\substack{ u=\left( u_{1},u_{2}\right)  \\ u_{1}\in
H^{1}\left( 0,2\pi \right) ,  \\ u_{2}\in L^{2}\left( 0,2\pi \right) }}\frac{%
\int_{0}^{2\pi }L\left( u\right) \cdot u\ ~dx}{\int_{0}^{2\pi }\rho _{0}%
\left[ u_{1}^{2}+u_{2}^{2}\right] ~dx}>0,  \label{lambda-k}
\end{equation}%
where 
\begin{eqnarray*}
L\left( u\right) \cdot u\ &=&-G\left( u_{1},u_{1x},u_{2}\right) +H\left(
u_{1}\right) \\
G\left( u_{1},u_{1x},u_{2}\right) &=&-\left( \gamma p_{0}+B_{0}^{2}\right)
\left( u_{1x}+ku_{2}+\frac{g\rho _{0}}{\gamma p_{0}+B_{0}^{2}}u_{1}\right)
^{2}, \\
H\left( u_{1}\right) &=&\left( \frac{g^{2}\rho _{0}^{2}}{\gamma
p_{0}+B_{0}^{2}}-g\rho _{0x}\right) u_{1}^{2},
\end{eqnarray*}%
The spectral radius $\Lambda $ is given by

\begin{equation}
\Lambda ^{2}=\sup_{v\in L^{2}\left( D\right) }\frac{\iint_{D}H\left(
v\right) dxdy}{\iint_{D}\rho _{0}v^{2}~dxdy}>0.  \label{Lambda2}
\end{equation}%
The sufficient condition for nonlinear instability is given by, for some $%
x_{0}\in \left( 0,2\pi \right) ,\,$

\begin{equation}
\frac{g^{2}\rho _{0}^{2}}{\gamma p_{0}+B_{0}^{2}}\left( x_{0}\right) >g\rho
_{0x}\left( x_{0}\right) .  \label{instability-criterion}
\end{equation}%
Notice that density inversion is not required for instability unlike the
incompressible case. Thus, even if light fluid is on top of heavy fluid and
magnetic field lines stabilizes, such steady state flows can't sustain
themselves against small initial disturbances under the effect of
compressibility. It means that the destabilizing effect of compressibility
dominates over the stabilizing one of magnetic field lines. Furthermore,
letting $\gamma \rightarrow \infty $ yields exactly (\ref{rho_0}) and the
instability criterion for the incompressible case can be recovered from
compressible ones as the limiting case.

(3) Magnetic field is parallel to the gravity $B_{0}\parallel \vec{g}$ for
incompressible fluids: 
\begin{equation}
B_{0}=\left( B_{0},0,0\right) .  \label{B2}
\end{equation}%
In this case we consider density profile $\rho _{0}$ whose gradient has
negative average over $\left( 0,2\pi \right) ,$ which is stronger than (\ref%
{rho_0}):%
\begin{equation*}
\int_{0}^{2\pi }-g\rho _{0x}\left( x\right) dx>0,\text{ }
\end{equation*}%
Linear operator $L$, Sturm-Liouville equation, and $\lambda _{k}$ are given
by:

\begin{equation}
L\left( v\right) :=B_{0}^{2}v_{xx}-\nabla \left( p_{t}+B_{0}\cdot Q\right) -%
\vec{g}\rho _{0x}v_{1},  \notag
\end{equation}

\begin{equation}
\rho _{0}\lambda _{k}^{2}\tilde{v}_{1}=B_{0}^{2}\tilde{v}_{1xx}+\frac{%
\lambda _{k}^{2}}{k^{2}}\left( \rho _{0}\tilde{v}_{1x}\right) _{x}-\frac{%
B_{0}^{2}}{k^{2}}\tilde{v}_{1xxxx}-g\rho _{0x}\tilde{v}_{1},  \label{SL1}
\end{equation}

\begin{equation}
\lambda _{k}^{2}=\sup_{u\in H^{2}\left( 0,2\pi \right) }\frac{\int_{0}^{2\pi
}L\left( u\right) \cdot u\ \ dx}{\int_{0}^{2\pi }\rho _{0}\left[ u^{2}+\frac{%
u_{x}^{2}}{k^{2}}\right] \,dx}>0,  \label{lambdak}
\end{equation}%
where

\begin{eqnarray*}
L\left( u\right) \cdot u\ &=&-F\left( -\frac{u_{xx}}{k}\right) +H\left(
u,u_{x}\right) , \\
F\left( -\frac{u_{xx}}{k}\right) &=&B_{0}^{2}\frac{u_{xx}^{2}}{k^{2}},\
H\left( u,u_{x}\right) =-g\rho _{0x}u^{2}-B_{0}^{2}u_{x}^{2},
\end{eqnarray*}%
$\Lambda $ has the following formula:

\begin{equation}
\Lambda ^{2}=\sup_{v\in H^{1}\left( D\right) }\frac{\iint_{D}H\left(
v,v_{x}\right) dxdy}{\iint_{D}\rho _{0}v^{2}\,dxdy}>0.  \label{Lambda1}
\end{equation}

\begin{lemma}
$\Lambda ^{2}$ in (\ref{Lambda1}) is positive.
\end{lemma}

\begin{proof}
We can take a family of test periodic functions in $H^{1}$ which guarantees
the positivity of $\Lambda $. Let $c=\min_{x}\rho _{0}\left( x\right) $ and 
\begin{equation*}
\phi ^{n}\left( x,y\right) =c-\frac{x}{n}\text{ \ for }0\leq x<2\pi .
\end{equation*}%
Since%
\begin{eqnarray*}
\iint_{D}-g\rho _{0x}\left( \phi ^{n}\right) ^{2}dxdy &\rightarrow
&c^{2}\iint_{D}-g\rho _{0x}dxdy>0\text{ \ as }n\rightarrow \infty \text{ },
\\
B_{0}^{2}\iint_{D}\left( \phi _{x}^{n}\right) ^{2}dxdy &=&\frac{4\pi
^{2}B_{0}^{2}}{n^{2}}\rightarrow 0\text{ \ as }n\rightarrow \infty \text{,}
\end{eqnarray*}%
we choose large $n$ to get positivity of the numerator of $\Lambda ^{2}$.
\end{proof}

(4) Magnetic field is parallel to the gravity $B_{0}\parallel \vec{g}$ for
compressible fluids: In this case, we obtain linear stability instead of
instability and we have the following:

\begin{equation*}
L(v)=B_{0}^{2}v_{xx}-\nabla \left( B_{0}\cdot Q\right) -\nabla \left(
p_{0x}v_{1}+\gamma p_{0}\nabla \cdot v\right) -\vec{g}\left( \rho
_{0x}v_{1}+\rho _{0}\nabla \cdot v\right) ,
\end{equation*}%
\begin{equation*}
\lambda _{k}^{2}=\sup_{\substack{ u=\left( u_{1},u_{2}\right)  \\ u_{1}\in
H^{1}\left( 0,2\pi \right) ,  \\ u_{2}\in H^{1}\left( 0,2\pi \right) }}\frac{%
\int_{0}^{2\pi }L\left( u\right) \cdot u\ ~dx}{\int_{0}^{2\pi }\rho _{0}%
\left[ u_{1}^{2}+u_{2}^{2}\right] ~dx}<0,
\end{equation*}%
where

\begin{eqnarray*}
L\left( u\right) \cdot u\ &=&-F\left( u_{2x}\right) -G\left(
u_{1},u_{1x},u_{2}\right) +H\left( u_{1},u_{1x}\right) , \\
F\left( u_{2x}\right) &=&B_{0}^{2}u_{2x}^{2}, \\
G\left( u_{1},u_{1x},u_{2}\right) &=&\left( \gamma p_{0}+B_{0}^{2}\right) 
\left[ \frac{\gamma p_{0}}{\gamma p_{0}+B_{0}^{2}}u_{1x}+\frac{g\rho
_{0}u_{1}}{\gamma p_{0}+B_{0}^{2}}+ku_{2}\right] ^{2}, \\
H\left( u_{1},u_{1x}\right) &=&\frac{1}{\gamma p_{0}+B_{0}^{2}}\left[ g\rho
_{0}u_{1}-B_{0}^{2}u_{1x}\right] ^{2}-g\rho
_{0x}u_{1}^{2}-B_{0}^{2}u_{1x}^{2} \\
&=&-\frac{B_{0}^{2}}{\left( \gamma p_{0}+B_{0}^{2}\right) \gamma p_{0}}%
\left( g\rho _{0}u_{1}+\gamma p_{0}u_{1x}\right) ^{2}<0.
\end{eqnarray*}%
and

\begin{equation}
\Lambda ^{2}=\sup_{v\in H^{1}\left( D\right) }\frac{\iint_{D}H\left(
v,v_{x}\right) \,dxdy}{\iint_{D}\rho _{0}v^{2}\,dxdy}<0.  \label{Lambda3}
\end{equation}

\section{Energy estimate and approximate solution}

In this section, we construct approximate solutions using a method
originated by Grenier in \cite{G} and we do energy estimates for the fully
nonlinear MHD system in both incompressible and compressible cases.

We first construct approximate solutions. In our construction, $\delta >0$
is an arbitrary small parameter, and $\theta $ is a small but fixed positive
constant (independent of $\delta $). We fix $k_{0}$ with $\lambda =\lambda
_{k_{0}}$ (dominant eigenvalue) so that 
\begin{equation*}
\Lambda <2\lambda .
\end{equation*}%
We define $T^{\delta }$ by%
\begin{equation}
\theta =\delta \exp \left( \lambda T^{\delta }\right) ,  \label{theta}
\end{equation}%
or equivalently,%
\begin{equation}
T^{\delta }=\frac{1}{\lambda }\ln \frac{\theta }{\delta }.  \label{T-delta}
\end{equation}%
We may write the full system (\ref{full-eqn}) in vector form for $w=\left(
\sigma \left( t,x,y\right) ,v\left( t,x,y\right) ,p\left( t,x,y\right)
,B\left( t,x,y\right) \right) $:%
\begin{equation*}
w+A^{1}\left( w\right) \partial _{x}w+A^{2}\left( w\right) \partial
_{y}w+L\left( w\right) =F\left( w\right) .
\end{equation*}%
An approximate solution $w^{a}\left( t,x,y\right) =\left( \sigma ^{a}\left(
t,x,y\right) ,v^{a}\left( t,x,y\right) ,p^{a}\left( t,x,y\right)
,B^{a}\left( t,x,y\right) \right) $ is of the form%
\begin{equation}
w^{a}\left( t,x,y\right) =\sum_{j=1}^{N}\delta ^{j}r_{j}\left( t,x,y\right) ,
\label{expansion}
\end{equation}%
such that%
\begin{equation}
w^{a}+A^{1}\left( w^{a}\right) \partial _{x}w^{a}+A^{2}\left( w^{a}\right)
\partial _{y}w^{a}+L\left( w^{a}\right) =F\left( w^{a}\right) +R_{N}^{a}.
\label{approximate}
\end{equation}%
We show existence of such approximate solutions in the following lemma. The
key point is that we can choose a dominant eigenvalue $\lambda $ with $%
\Lambda <2\lambda $ to make this construction work.

\begin{lemma}
For any fixed integer $N>0$, there exists an approximate solution (\ref%
{expansion}) satisfying (\ref{approximate}). Furthermore, for every integer $%
s\geq 0$, there is small $\theta >0$ such that if $0\leq t\leq T^{\delta }$
as in (\ref{T-delta}), the $j$-th coefficient $r_{j}$ and the remainder $%
R_{N}^{a}$ satisfy%
\begin{align}
\left\Vert r_{j}\right\Vert _{H^{s}}& \leq C\left( s,N\right) \exp \left(
j\lambda t\right) ,~~\text{for }1\leq j\leq N,  \label{r-j} \\
\left\Vert R_{N}^{a}\right\Vert _{H^{s}}& \leq C\left( s,N\right) \delta
^{N+1}\exp \left\{ \left( N+1\right) \lambda t\right\} .  \label{R}
\end{align}
\end{lemma}

\begin{proof}
The construction of $r_{j}$ will be made by induction on $j$. The idea is as
follows. We split the system into linear and nonlinear part:%
\begin{align}
& \partial _{t}w^{a}+A^{1}\left( 0\right) \partial _{x}w^{a}+A^{2}\left(
0\right) \partial _{y}w^{a}+L\left( w^{a}\right)  \label{r-j-construction} \\
& =\left[ A^{1}\left( w^{a}\right) -A^{1}\left( 0\right) \right] \partial
_{x}w^{a}+\left[ A^{2}\left( u^{a}\right) -A^{2}\left( 0\right) \right]
\partial _{y}w^{a}-F\left( w^{a}\right) =:h\left( \delta \right)  \notag
\end{align}
with the Taylor expansion of $w^{a}$ in $\delta$%
\begin{equation*}
w^{a}=\sum_{j=1}^{N}\delta^{j}r_{j}.
\end{equation*}
Thus, our $r_{j}$ is the solution of the part of (\ref{r-j-construction})
which corresponds to the coefficient of $\delta^{j}$ in its Taylor expansion.

For $j=1$, take for $r_{1}$ the smooth normal growing mode to the linearized
system with our chosen wave number $k_{0}$ and the corresponding dominant
eigenvalue $\lambda=\lambda_{k_{0}}$ as in (\ref{expo}). Clearly this
growing mode fulfills (\ref{r-j}).

Assuming that we have constructed $r_{j}$ $(j<N)$ which satisfies (\ref{r-j}%
), we construct $r_{j+1}$. Let 
\begin{equation*}
u_{j}=\sum_{k=1}^{j}\delta ^{k}r_{k}\left( t,x,y\right) .
\end{equation*}%
Define the nonlinear part of the system substituted by $u_{j}$ as 
\begin{equation*}
h_{j+1}\left( \delta \right) =\left[ A^{1}\left( u_{j}\right) -A^{1}\left(
0\right) \right] \partial _{x}u_{j}+\left[ A^{2}\left( u_{j}\right)
-A^{2}\left( 0\right) \right] \partial _{y}u_{j}-F\left( u_{j}\right) .
\end{equation*}%
Since this is the nonlinear part of the system and the terms in $\delta
^{j+1}$ come from the terms $\delta ^{k}$ for $1\leq k\leq j$, it is enough
to consider $u_{j}$ in order to collect the $(j+1)$-th coefficient of
nonlinear part of the expansion. Then $r_{j+1}$ is defined to be a solution
of%
\begin{align*}
& \partial _{t}r_{j+1}+A^{1}\left( 0\right) \partial _{x}r_{j+1}+A^{2}\left(
0\right) \partial _{y}r_{j+1}+L\left( r_{j+1}\right) \\
& =\frac{-h_{j+1}^{\left( j+1\right) }\left( 0\right) }{\left( j+1\right) !}
\end{align*}%
with initial data $r_{j+1}\left( 0,x,y\right) =0$. Notice that%
\begin{equation*}
\frac{-h_{j+1}^{\left( j+1\right) }\left( 0\right) }{\left( j+1\right) !}%
=\sum_{j_{1}m_{1}+j_{2}m_{2}+..+j_{p+1}m_{p+1}=j+1}B_{J,i}^{M}\
r_{j_{1}}^{m_{1}}r_{j_{2}}^{m_{2}}\cdot \cdot r_{j_{p}}^{m_{p}}\partial
_{i}r_{j_{p+1}}^{m_{p+1}},
\end{equation*}%
where $m_{k}\geq 0,\ 1\leq k\leq j,$ and $B_{J,i}^{M}$ depends on $A^{i}$
and $F$. Induction hypothesis (\ref{r-j}) for $r_{k},\ 1\leq k\leq j$
applies to get, for all $s$,%
\begin{equation*}
\left\Vert \frac{h_{j+1}^{\left( j+1\right) }\left( 0\right) }{\left(
j+1\right) !}\right\Vert _{H^{s}}\leq C\left( s,N\right) \exp \left[ \left(
j_{1}m_{1}+j_{2}m_{2}+..+j_{p+1}m_{p+1}\right) \lambda t\right] =C\left(
s,N\right) \exp \left[ \left( j+1\right) \lambda t\right] .
\end{equation*}%
Thanks to our linear estimates for $\Lambda $ in Lemma 4 and Duhamel
principle, we have%
\begin{align*}
\left\Vert r_{j+1}\left( t,\cdot \right) \right\Vert _{H^{s}}& \leq
C\int_{0}^{t}e^{\Lambda \left( t-\tau \right) }\left\Vert \frac{%
h_{j+1}^{\left( j+1\right) }\left( 0\right) }{\left( j+1\right) !}\left(
\tau \right) \right\Vert _{H^{s+2}} \\
& \leq C\left( s,N\right) \int_{0}^{t}e^{\Lambda \left( t-\tau \right)
}e^{\left( j+1\right) \lambda \tau }d\tau \\
& \leq C\left( s,N\right) e^{\left( j+1\right) \lambda t}
\end{align*}%
since $j+1\geq 2$ and $\Lambda <2\lambda $.

We now define $w^{a}=\sum_{j=1}^{N}\delta^{j}r_{j}$ and it satisfies%
\begin{equation*}
\partial_{t}w^{a}+A^{1}\left( 0\right) \partial_{x}w^{a}+A^{2}\partial
_{y}w^{a}+L\left( w^{a}\right) =-\sum_{j=1}^{N}\frac{\delta^{j}h_{j}^{\left(
j\right) }\left( 0\right) }{j!}.
\end{equation*}
Then $R_{N}^{a}$ is defined to be the sum of all higher terms than $N$ in
nonlinear part of the $\delta$-expansion (\ref{r-j-construction}):%
\begin{equation*}
R_{N}^{a}=h\left( \delta\right) -\sum_{j=1}^{N}\frac{\delta^{j}h_{j}^{\left(
j\right) }\left( 0\right) }{j!},
\end{equation*}
which clearly satisfies (\ref{R}) and our proof is complete.
\end{proof}

We state local in time existence for the incompressible ideal MHD equations:

\begin{lemma}
(Local existence to the full system) For all $s\geq3$ and for any given
initial data $\left( \sigma_{0},v_{0}\right) \in H^{s}\left( D\right) $ such
that $\rho\left( 0\right) \equiv\rho_{0}\left( x\right) +\sigma _{0}\left(
x,y\right) \geq m>0,$ there is a $T>0$ such that there exists a unique
solution $\left( \sigma,v,p,B\right) \in C\left( \left[ 0,T\right]
:H^{s}\left( D\right) \right) $ to (\ref{full-eqn}) with $\rho\left(
t\right) =\rho_{0}\left( x\right) +\sigma\left( t,x,y\right) >0.$
\end{lemma}

We first treat the incompressible case and we use some vector identities to
get:%
\begin{equation*}
\sigma _{t}+v\cdot \nabla \left( \rho _{0}+\sigma \right) =0,
\end{equation*}%
\begin{align*}
\left( \rho _{0}+\sigma \right) \left( v_{t}+v\cdot \nabla v\right) &
=-\nabla \left( p+\frac{1}{2}B^{2}+B_{0}\cdot B\right) +B\cdot \nabla B \\
& +B_{0}\cdot \nabla B+B\cdot \nabla B_{0}+\sigma \vec{g},
\end{align*}%
\begin{equation*}
B_{t}=B_{0}\cdot \nabla v-v\cdot \nabla B_{0}+B\cdot \nabla v-v\cdot \nabla
B,
\end{equation*}%
where $\hat{k}$ is the unit vector in the $z$ direction. Let $w\left(
t,x,y\right) =\left( \sigma \left( t,x,y\right) ,v\left( t,x,y\right)
,p\left( t,x,y\right) ,B\left( t,x,y\right) \right) \in C\left( \left[ 0,T%
\right] ;H^{s}\left( D\right) \right) $ be a local solution as constructed
above.

Let $w^{a}\left( t,x,y\right) =\left( \sigma^{a}\left( t,x,y\right)
,v^{a}\left( t,x,y\right) ,p^{a}\left( t,x,y\right) ,B^{a}\left(
t,x,y\right) \right) $ be an approximate solution. We now estimate the
difference 
\begin{equation*}
w^{d}\left( t,x,y\right) =w\left( t,x,y\right) -w^{a}\left( t,x,y\right) .
\end{equation*}%
\begin{equation}
\sigma_{t}^{d}+\nabla\sigma^{d}\cdot v+\nabla\left( \sigma^{a}+\rho
_{0}\right) \cdot v^{d}=-R_{N,1}^{a},  \label{diff1}
\end{equation}%
\begin{align}
& \left( \rho_{0}+\sigma\right) \left( v_{t}^{d}+v^{d}\cdot\nabla
v+v^{a}\cdot\nabla v^{d}\right) +\sigma^{d}\left( v_{t}^{a}+v^{a}\cdot\nabla
v^{a}\right)  \notag \\
& =-\nabla\left( p^{d}+\frac{1}{2}B^{d}\left( B+B^{a}\right) +B_{0}\cdot
B^{d}\right) +B^{d}\cdot\nabla B+B^{a}\cdot\nabla B^{d}  \label{diff2} \\
& +B_{0}\cdot\nabla B^{d}+B^{d}\cdot\nabla B_{0}+\sigma^{d}\vec{g}%
-R_{N,2}^{a},  \notag
\end{align}

\begin{equation}
B_{t}^{d}=B_{0}\cdot\nabla v^{d}-v^{d}\cdot\nabla B_{0}+B^{d}\cdot\nabla
v+B^{a}\cdot\nabla v^{d}-v^{d}\cdot\nabla B-v^{a}\cdot\nabla
B^{d}-R_{N,3}^{a}.  \label{diff3}
\end{equation}

\begin{lemma}
For $s\geq3$ and assume $\left\Vert \sigma\right\Vert _{\infty}\leq\frac{1}{2%
}\left\Vert \rho_{0}\right\Vert _{\infty}$, then there exists a continuous
positive function $C$ depending only on $s,\rho_{0}$ such that%
\begin{align*}
\frac{d}{dt}\left\Vert w^{d}\right\Vert _{s}^{2} & \leq C\left( \left\Vert
w^{d}\right\Vert _{s},\left\Vert w^{a}\right\Vert _{s+1}\right) \left\Vert
w^{d}\right\Vert _{s}^{2} \\
& +\left\Vert R_{N,1}^{a}\right\Vert _{s}^{2}+\left\Vert
R_{N,2}^{a}\right\Vert _{s}^{2}+\left\Vert R_{N,3}^{a}\right\Vert _{s}^{2}.
\end{align*}
\end{lemma}

\begin{proof}
This energy estimate is straightforward and thus we give a brief sketch. We
take $\partial _{\alpha }$ of the equations (\ref{diff1})-(\ref{diff3}),
multiply through $\partial _{\alpha }$-derivatives, and integrate over $D$.
For symmetric terms such as $\nabla \partial _{\alpha }\sigma ^{d}\cdot v$
in the first equation vanish upon integration due to the divergence free
condition for $v,v^{a},v^{d},B,B^{a},B^{d}$. We now estimate nonsymmetric
terms. Terms which need attention are $\left( B^{d}\cdot \nabla \right)
\partial _{\alpha }B^{d},\left( B^{a}\cdot \nabla \right) \partial _{\alpha
}B^{d},\left( B_{0}\cdot \nabla \right) \partial _{\alpha }B^{d}$ from $v$%
-equation (\ref{diff2}) and $\left( B^{d}\cdot \nabla \right) \partial
_{\alpha }v^{d},\left( B^{a}\cdot \nabla \right) \partial _{\alpha }v^{d},$ $%
\left( B_{0}\cdot \nabla \right) \partial _{\alpha }B^{d}$ from $B$-equation
(\ref{diff3}). Upon integration, the three corresponding pairs exactly
cancel out. To see this,%
\begin{align*}
B_{i}^{d}\partial _{i}\partial _{\alpha }B_{j}^{d}\partial _{\alpha
}v_{j}^{d}& =\partial _{i}B_{i}^{d}\partial _{\alpha }B_{j}^{d}\partial
_{\alpha }v_{j}^{d}-B_{i}^{d}\partial _{i}\partial _{\alpha
}v_{j}^{d}\partial _{\alpha }B_{j}^{d}=-B_{i}^{d}\partial _{i}\partial
_{\alpha }v_{j}^{d}\partial _{\alpha }B_{j}^{d}, \\
B_{i}^{a}\partial _{i}\partial _{\alpha }B_{j}^{d}\partial _{\alpha
}v_{j}^{d}& =\partial _{i}B_{i}^{a}\partial _{\alpha }B_{j}^{d}\partial
_{\alpha }v_{j}^{d}-B_{i}^{a}\partial _{i}\partial _{\alpha
}v_{j}^{d}\partial _{\alpha }B_{j}^{d}=-B_{i}^{a}\partial _{i}\partial
_{\alpha }v_{j}^{d}\partial _{\alpha }B_{j}^{d}, \\
B_{0i}\partial _{i}\partial _{\alpha }B_{j}^{d}\partial _{\alpha }v_{j}^{d}&
=\partial _{i}B_{0i}\partial _{\alpha }B_{j}^{d}\partial _{\alpha
}v_{j}^{d}-B_{0i}\partial _{i}\partial _{\alpha }v_{j}^{d}\partial _{\alpha
}B_{j}^{d}=-B_{0i}\partial _{i}\partial _{\alpha }v_{j}^{d}\partial _{\alpha
}B_{j}^{d},
\end{align*}%
thanks to divergence free condition for $v$ and $B$. Therefore we obtain the
Lemma.
\end{proof}

We now extend the energy estimates to the compressible case using
symmetrizer, which is necessary in this compressible case. We write the full
system (\ref{full-eqn}) in components: 
\begin{equation*}
\sigma _{,t}+v_{i}\sigma _{,i}+\left( \rho _{0}+\sigma \right) v_{i,i}+\rho
_{0x}v_{1}=0,
\end{equation*}%
\begin{eqnarray*}
&&\left( \rho _{0}+\sigma \right) \left\{ v_{j,t}+v_{i}v_{j,i}\right\}
+\left( \rho _{0}+\sigma \right) q\sigma
_{,j}+B_{i}B_{i,j}-B_{i}B_{j,i}+B_{0}B_{3,j} \\
&=&\left[ \rho _{0}q\left( \rho _{0}\right) -\left( \rho _{0}+\sigma \right)
q\right] \rho _{0,j}+B_{i}B_{0,i}-B_{0,j}B_{3}
\end{eqnarray*}

\begin{equation*}
B_{j,t}-B_{i}v_{j,i}+v_{i,i}B_{j}+v_{i}B_{j,i}+v_{i,i}B_{0}=-v_{i}B_{0,i}.
\end{equation*}%
Here $j=1,2,3$ and twice $i$ means the sum over $i=1,2,3$.

We rewrite the full system near the steady state $\left( \rho _{0},\vec{0}%
,B\right) $ in vector notations for $w=\left( \sigma
,v_{1},v_{2},v_{3},B_{1},B_{2},B_{3}\right) $:

\begin{equation}
w_{t}+A^{1}\left( w\right) \partial _{x}w+A^{2}\left( w\right) \partial
_{y}w+L\left( w\right) =F\left( w\right) ,  \label{u}
\end{equation}%
Introducing and multiplying the symmetrizer with $q\left( \rho \right)
=\gamma \rho ^{\gamma -2}$ 
\begin{equation}
D=\text{diag}\left( q\left( \rho _{0}+\sigma \right) ,\rho _{0}+\sigma ,\rho
_{0}+\sigma ,\rho _{0}+\sigma ,1,1,1\right) ,  \label{symmetrizer}
\end{equation}%
leads to the following symmetric matrices%
\begin{equation*}
DA_{1}=\left[ 
\begin{array}{ccccccc}
q\left( \rho _{0}+\sigma \right) v_{1} & \left( \rho _{0}+\sigma \right)
q\left( \rho _{0}+\sigma \right) & 0 & 0 & 0 & 0 & 0 \\ 
\left( \rho _{0}+\sigma \right) q\left( \rho _{0}+\sigma \right) & \left(
\rho _{0}+\sigma \right) v_{1} & 0 & 0 & 0 & B_{2} & B_{3}+B_{0} \\ 
0 & 0 & \left( \rho _{0}+\sigma \right) v_{1} & 0 & 0 & -B_{1} & 0 \\ 
0 & 0 & 0 & \left( \rho _{0}+\sigma \right) v_{1} & 0 & 0 & -B_{1} \\ 
0 & 0 & 0 & 0 & v_{1} & 0 & 0 \\ 
0 & B_{2} & -B_{1} & 0 & 0 & v_{1} & 0 \\ 
0 & B_{3}+B_{0} & 0 & -B_{1} & 0 & 0 & v_{1}%
\end{array}%
\right] ,
\end{equation*}%
\begin{equation*}
DA_{2}=\left[ 
\begin{array}{ccccccc}
q\left( \rho _{0}+\sigma \right) v_{2} & 0 & \left( \rho _{0}+\sigma \right)
q\left( \rho _{0}+\sigma \right) & 0 & 0 & 0 & 0 \\ 
0 & \left( \rho _{0}+\sigma \right) v_{2} & 0 & 0 & -B_{2} & 0 & 0 \\ 
\left( \rho _{0}+\sigma \right) q\left( \rho _{0}+\sigma \right) & 0 & 
\left( \rho _{0}+\sigma \right) v_{2} & 0 & B_{1} & 0 & B_{3}+B_{0} \\ 
0 & 0 & 0 & \left( \rho _{0}+\sigma \right) v_{2} & 0 & 0 & -B_{2} \\ 
0 & -B_{2} & B_{1} & 0 & v_{2} & 0 & 0 \\ 
0 & 0 & 0 & 0 & 0 & v_{2} & 0 \\ 
0 & 0 & B_{3}+B_{0} & -B_{2} & 0 & 0 & v_{2}%
\end{array}%
\right] ,
\end{equation*}

\begin{equation*}
DA_{3}=\left[ 
\begin{array}{ccccccc}
q\left( \rho _{0}+\sigma \right) v_{3} & 0 & 0 & \left( \rho _{0}+\sigma
\right) q\left( \rho _{0}+\sigma \right) & 0 & 0 & 0 \\ 
0 & \left( \rho _{0}+\sigma \right) v_{3} & 0 & 0 & -B_{3} & 0 & 0 \\ 
0 & 0 & \left( \rho _{0}+\sigma \right) v_{3} & 0 & 0 & -B_{3} & 0 \\ 
\left( \rho _{0}+\sigma \right) q\left( \rho _{0}+\sigma \right) & 0 & 0 & 
\left( \rho _{0}+\sigma \right) v_{3} & B_{1} & B_{2} & B_{0} \\ 
0 & -B_{3} & 0 & B_{1} & v_{3} & 0 & 0 \\ 
0 & 0 & -B_{3} & B_{2} & 0 & v_{3} & 0 \\ 
0 & 0 & 0 & B_{0} & 0 & 0 & v_{3}%
\end{array}%
\right] ,
\end{equation*}

\begin{equation*}
DL\left( u\right) =\left[ 
\begin{array}{c}
\rho _{0x}q\left( \rho _{0}\right) v_{1} \\ 
q\left( \rho _{0}\right) \rho _{0x}\sigma +\rho _{0}q^{\prime }\left( \rho
_{0}\right) \rho _{0x}\sigma +B_{0x}B_{3} \\ 
0 \\ 
-B_{0x}B_{1} \\ 
0 \\ 
0 \\ 
B_{0x}v_{1}%
\end{array}%
\right] \ ,
\end{equation*}

\begin{equation*}
\ DF\left( u\right) =\left[ 
\begin{array}{c}
\rho _{0x}q\left( \rho _{0}\right) v_{1}-\rho _{0x}q\left( \rho _{0}+\sigma
\right) v_{1} \\ 
\left[ \left( \sigma +\rho _{0}\right) q\left( \sigma +\rho _{0}\right)
-\rho _{0}q\left( \rho _{0}\right) -q\left( \rho _{0}\right) \rho
_{0x}\sigma -\rho _{0}q^{\prime }\left( \rho _{0}\right) \rho _{0x}\sigma %
\right] \rho _{0x} \\ 
0 \\ 
0 \\ 
0 \\ 
0 \\ 
0%
\end{array}%
\right] .
\end{equation*}%
Now that $DA^{1},$ $DA^{2}$ and $DA^{3}$ are symmetric, we have the
following local in time solution to the full system via standard energy
estimates for symmtrizable hyperbolic system as in \cite{CGG}, \cite{M} in
the absence of magnetic fields.

\begin{lemma}
For all $s\geq 3$ and for any given initial data $\left( \sigma
_{0},v_{0}\right) \in H^{s}\left( D\right) $ such that $\rho \left( 0\right)
\equiv \rho _{0}\left( x\right) +\sigma _{0}\left( x,y\right) \geq m>0$,
there is a $T>0$ such that there exists a unique solution $u=\left( \sigma
,v,B\right) \in C\left( \left[ 0,T\right] :H^{s}\left( D\right) \right) $ to
(\ref{full-eqn}) with $\rho \left( t\right) =\rho _{0}\left( x\right)
+\sigma \left( t,x,y\right) >0$.
\end{lemma}

Since construction of approximate solutions is similar to the incompressible
case as in Lemma 5 for a hyperbolic system, we omit it. We now estimate the
difference of \ an exact solution and an approximate solution. Let $w\left(
t,x,y\right) =\left( \sigma \left( t,x,y\right) ,v\left( t,x,y\right)
,B\left( t,x,y\right) \right) \in C\left( \left[ 0,T\right] :H^{s}\left(
D\right) \right) $ be an exact solution and $w^{a}\left( t,x,y\right)
=\left( \sigma ^{a}\left( t,x,y\right) ,v^{a}\left( t,x,y\right)
,B^{a}\left( t,x,y\right) \right) $ be an approximate solution as
constructed in Lemma 5. Then their difference is%
\begin{equation*}
w^{d}=w-w^{a}=\left( \sigma -\sigma ^{a},v-v^{a},B-B^{a}\right) ,
\end{equation*}%
and it satisfies%
\begin{eqnarray}
&&w_{t}^{d}+\sum_{i=1}^{3}A^{i}\left( w^{a}+w^{d}\right) \partial
_{i}w^{d}+\sum_{i=1}^{3}\left[ A^{1}\left( w^{a}+w^{d}\right) -A^{1}\left(
w^{a}\right) \right] \partial _{x}w^{a}  \label{difference} \\
&=&-L\left( w^{d}\right) +F\left( w\right) -F\left( w^{a}\right) -R_{N}^{a}.
\notag
\end{eqnarray}%
This symmetrizable hyperbolic system for $w^{d}$ allows the following energy
estimates. The proof is straightforward by classical energy methods as in 
\cite{CGG}, \cite{G}.

\begin{lemma}
For any $s\geq 3$, there exists a continuous function $g_{s}\left( \cdot
,\cdot \right) $ such that%
\begin{equation*}
\frac{\partial }{\partial t}\left\vert \left\vert \left\vert
w^{d}\right\vert \right\vert \right\vert _{s}^{2}\leq g_{s}\left( \left\vert
\left\vert \left\vert w^{d}\right\vert \right\vert \right\vert
_{s},\left\vert \left\vert \left\vert w^{a}\right\vert \right\vert
\right\vert _{s+1}\right) \left\vert \left\vert \left\vert w^{d}\right\vert
\right\vert \right\vert _{s}^{2}+\left\vert \left\vert \left\vert
R_{N}^{a}\right\vert \right\vert \right\vert _{s}^{2},
\end{equation*}

where $\left\vert \left\vert \left\vert \cdot \right\vert \right\vert
\right\vert _{s}$ is defined by%
\begin{equation*}
\left\vert \left\vert \left\vert v\right\vert \right\vert \right\vert
_{s}^{2}=\sum_{\left\vert \alpha \right\vert \leq s}\partial _{\alpha
}vD\left( \rho _{0}+\sigma \right) \partial _{\alpha }v.
\end{equation*}
\end{lemma}

Here $D$ is the symmetrizer as in (\ref{symmetrizer}) and notice that $%
\left\vert \left\vert \left\vert \cdot \right\vert \right\vert \right\vert
_{s}$ is related to the usual norm $\left\vert \left\vert \cdot \right\vert
\right\vert _{H^{s}}$ by 
\begin{equation*}
\eta \left\Vert v\right\Vert _{H^{s}}\leq \left\vert \left\vert \left\vert
v\right\vert \right\vert \right\vert _{s}\leq C_{s}\left( \left\Vert \sigma
\right\Vert _{H^{s}},\left\Vert \rho _{0}\right\Vert _{H^{s}}\right)
\left\Vert v\right\Vert _{H^{s}},
\end{equation*}

since $\rho _{0}\geq c>0$ and so $D\geq \eta Id$ for some $\eta $ and $s\geq
3$.

Notice that all three norms $\left\Vert \cdot \right\Vert _{s}$, $\left\vert
\left\vert \left\vert \cdot \right\vert \right\vert \right\vert _{s},$ and $%
\left\Vert \cdot \right\Vert _{H^{s}}$ are equivalent since $\rho _{0}$ is
smooth with a positive minimum.

\section{Nonlinear instability}

\begin{proof}[Proof of Theorem 2]
Let $w^{a}\left( t,x,y\right) =\left( \sigma ^{a}\left( t,x,y\right)
,v^{a}\left( t,x,y\right) ,p^{a}\left( t,x,y\right) ,B^{a}\left(
t,x,y\right) \right) $ be an approximate solution with $N$ to be determined
later. For any $\delta >0$, there exists a local in time solution $w^{\delta
}\left( t,x,y\right) =$ $\left( \sigma ^{\delta }\left( t,x,y\right)
,v^{\delta }\left( t,x,y\right) ,p^{\delta }\left( t,x,y\right) ,B^{\delta
}\left( t,x,y\right) \right) $ with initial data $w^{a}\left( 0\right) $ to
the full system (\ref{full-eqn}). By the Lemma 7 and Lemma 9, we have%
\begin{equation}
\frac{d}{dt}\left\Vert w^{d}\right\Vert _{s}^{2}\leq C\left( \left\Vert
w^{d}\right\Vert _{s},\left\Vert w^{a}\right\Vert _{s+1}\right) \left\Vert
w^{d}\right\Vert _{s}^{2}+C\delta ^{2\left( N+1\right) }e^{2\left(
N+1\right) \lambda t} ,  \label{wd}
\end{equation}%
with $w^{d}\left( 0\right) =0$.

Let%
\begin{equation*}
T=\sup\left\{ t\ |\ \left\Vert w^{a}\right\Vert _{s+1}\leq\omega /2,\
\left\Vert w^{d}\right\Vert _{s}\leq\omega/2\right\} ,
\end{equation*}
where $\omega$ is a small positive number which assures local existence.
Then $T$ depending on $\delta$ is well defined since $w^{d}\left( 0\right)
=0 $ and $\left\Vert w^{\alpha}\left( 0\right) \right\Vert _{s+1}=O\left(
\delta\right) $. We claim that the instability time $T^{\delta}$ occurs
within the existence time $T$, that is, $T^{\delta}\leq T$. Suppose not,
i.e., $T<T^{\delta}$. Then for $t\leq T$, by our construction of approximate
solution, we have%
\begin{equation*}
\left\Vert w^{a}\right\Vert _{s+1}\leq C\sum_{j=1}^{N}\delta^{j}\left\Vert
r_{j}\left( t\right) \right\Vert _{H^{s}}\leq\sum_{j=1}^{N}C_{j}\delta
^{j}e^{j\lambda t}\leq\sum_{j=1}^{N}C_{j}\theta^{j}<\omega/2\text{ },
\end{equation*}
if $\theta$ is small. Now we appeal to the definition of $T$ and (\ref{wd})
to get, for $t\leq T$,%
\begin{equation*}
\frac{d}{dt}\left\Vert w^{d}\right\Vert _{s}^{2}\leq C\left( \omega
/2,\omega/2\right) \left\Vert w^{d}\right\Vert _{s}^{2}+C\delta^{2\left(
N+1\right) }e^{2\left( N+1\right) \lambda t}.
\end{equation*}
We choose $N>0$ satisfying%
\begin{equation*}
C\left( \omega/2,\omega/2\right) <2\left( N+1\right) ,
\end{equation*}
so that $\left\Vert w^{d}\right\Vert _{s}$ has growth rate at most $\left(
N+1\right) \lambda$. Then using Gronwall inequality leads to%
\begin{equation*}
\left\Vert w^{d}\right\Vert _{s}\leq C\delta^{2\left( N+1\right) }e^{2\left(
N+1\right) \lambda t}=C\theta^{2\left( N+1\right) }<\omega/2,
\end{equation*}
if $\omega $ is small. Thus we can deduce $T^{\delta }\leq T$. Now at the
instability time $T^{\delta }$,%
\begin{align*}
\left\Vert w^{a}\left( T^{\delta }\right) \right\Vert _{L^{1}}& \geq \delta
\left\Vert r_{1}\left( T^{\delta }\right) \right\Vert
_{L^{1}}-\sum_{j=2}^{N}\delta ^{j}\left\Vert r_{j}\left( T^{\delta }\right)
\right\Vert _{L^{1}} \\
& \geq C\delta e^{\lambda T^{\delta }}-\sum_{j=2}^{N}C_{j}\delta
^{j}e^{j\lambda T^{\delta }} \\
& =C\theta -\sum_{j=2}^{N}C_{j}\theta ^{j}\geq \frac{C}{2}\theta .
\end{align*}%
We then deduce at time $T^{\delta }$,%
\begin{align*}
\left\Vert w^{\delta }\left( T^{\delta }\right) \right\Vert _{L^{1}}& \geq
\left\Vert w^{\delta }\left( T^{\delta }\right) \right\Vert
_{L^{1}}-\left\Vert \left( w^{\delta }-w^{a}\right) \left( T^{\delta
}\right) \right\Vert _{L^{1}} \\
& \geq \left\Vert w^{a}\left( T^{\delta }\right) \right\Vert
_{L^{1}}-\left\Vert \left( w^{\delta }-w^{a}\right) \left( T^{\delta
}\right) \right\Vert _{H^{s}} \\
& \geq \frac{C}{2}\theta -C\theta ^{N+1}\geq \frac{C}{4}\theta =\varepsilon
_{0}>0.
\end{align*}
\end{proof}

\end{document}